\documentclass[a4paper,12pt]{article}

\usepackage{latexsym}
\usepackage{amssymb}
\usepackage{amsthm}
\usepackage{amscd}
\usepackage{amsmath}
\usepackage{graphics}

\newtheorem{theorem}{Theorem}
\newtheorem{lemma}{Lemma}
\newtheorem{corollary}{Corollary}
\newtheorem{proposition}{Proposition}
\newtheorem{remark}{Remark}

\newtheorem{example}{Example}

\newcommand{\fkp}{\mathfrak p}

\newcommand{\fkm}{\mathfrak m}

\newcommand{\fkn}{\mathfrak n}

\DeclareMathOperator{\Spec}{\mathrm{Spec}}
\DeclareMathOperator{\Supp}{\mathrm{Supp}}

\DeclareMathOperator{\Proj}{\mathrm{Proj}}
\DeclareMathOperator{\PProj}{\mathbf{Proj}}
\newcommand{\Exc}{\mathrm{Excess}}

\newcommand{\Rees}{\mathcal R}
\newcommand{\Oh}{\mathcal O}
\newcommand{\adj}{\mathcal J}
\newcommand{\R}{R}
\newcommand{\B}{\mathbf B}

\newcommand{\ord}{\mathrm{ord}}

\begin{document}

\title{Adjoint ideals and Gorenstein blowups in two-dimensional 
regular local rings}

\author{Eero Hyry\
\footnote{
University of Helsinki, Department of Mathematics and Statistics, P.~O.~Box 68, FIN--00014 University of Helsinki, Finland,
eero.hyry@helsinki.fi}
 \and 
Yukio Nakamura
\footnote{
Department of Mathematics, School of Science and Technology, Meiji University, 214-8571, Japan,
ynakamu@math.meiji.ac.jp}
\and 
Lauri Ojala
\footnote{
University of Helsinki, Department of Mathematics and Statistics, P.~O.~Box 68, FIN--00014 University of Helsinki, Finland,
lauri.ojala@helsinki.fi}}

\maketitle

\thanks{Dedicated to Prof.~Kei-ichi Watanabe on the occasion of his
60th birthday.}

\begin{abstract}

In this article we investigate when a complete ideal in a two-dimensional regular local
ring is a multiplier ideal of some ideal with an integral multiplying parameter. In particular,
we show that this question is closely connected to the Gorenstein property of the blowup
along the ideal.

\end{abstract}

\section{Introduction}
\label{sec:intro}

Multiplier ideals have recently emerged as a fundamental tool in algebraic 
geometry.
To a given ideal $I$ in a regular local ring one can attach a collection of 
multiplier ideals 
$\adj(cI)$ depending on a rational multiplying parameter $c$ 
(see \cite{Laz}).  Multiplier ideals 
corresponding to an integer value of the parameter were independently 
introduced into commutative algebra
by Lipman who called them adjoint ideals (see \cite{L3} and \cite{L5}).  
As multiplier ideals are always integrally 
closed, it is natural to ask how special multiplier ideals are among all 
integrally closed ideals. In dimension two this 
question was answered by Lipman and Watanabe in \cite{LW}, and independently 
by Favre and Jonsson in \cite{FJ}. 
They showed that every integrally closed (complete) ideal in a 
two-dimensional regular local ring is a multiplier ideal of 
some ideal for some value of the multiplying parameter. However, 
even in this case it remains open in general when this
parameter can be taken to be an integer i.e.~the ideal is 
an adjoint ideal of some ideal. 

The purpose of this article is to investigate when a complete ideal 
$J$ in a two-dimensional local ring $(A,\fkm)$ is an adjoint ideal.  
This question can be rephrased in more algebraic terms by using a  
result of Huneke and Swanson saying that the adjoint 
$\adj(I)$ of any complete $\fkm$-primary ideal $I$ coincides with the second 
Fitting-ideal $F_2(I)$ (see~\cite[Proposition 3.16]{HuSw}). One can then equivalently ask which complete ideals 
$J$ are second Fitting-ideals.  
Our arguments are based on the Zariski-Lipman theory of complete ideals.  
In particular, they rely strongly on the notion of proximity. 
In this framework, one can define the adjoint ideal $\adj(I)$ of a complete 
$\fkm$-primary ideal $I\subset A$ as the unique complete ideal whose order at an infinitely 
near point is one less than that of $I$, 
if the latter is greater than one, and zero elsewhere. 

We will now describe our main results in more detail. We first observe in
Lemma~\ref{firstlem} that it is enough to consider $\fkm$-primary ideals. In Theorem \ref{main} 
we then show that a complete
$\fkm$-primary ideal $J$ is an adjoint ideal if and only if 
$\adj(J)$ divides $J^2$ i.e.~$F_2(J)$ divides $J^2$. If this is
the case, then higher powers of $J$ turn out to be iterated 
adjoints. In fact, by a further result of Huneke and Swanson they
are then higher Fitting ideals (see Proposition~\ref{Fittingprop}).
Our second main observation connects the property of being an adjoint ideal 
to the Gorenstein property of the blowup along the ideal. 
Our Theorem \ref{gor} says that $J$ 
is an adjoint if and only if the blowup $Y=\Proj R_A(J)$ along $J$ is 
Gorenstein and the sheaf
$J\omega_Y^{-1}$ where $\omega_Y$ denotes the canonical sheaf of 
$Y$, is globally generated.  In particular, this implies that $\Proj R_A(J)$ is Gorenstein 
if and only if $J^n$ is an adjoint ideal for large $n$. Actually, it is enough to
check that $J^{\ord_A(J)-1}$ is an adjoint ideal or equivalently $\adj(J)$ divides
$J^{\ord_A(J)}$, where $\ord_A(J)$ denotes the order of $J$ (see Corollary \ref{gor_ord}).

Significantly, $\Proj R_A(\adj(J))$ is always 
Gorenstein as well as 
$\Proj R_A(\adj(J^2))$. Since $\adj(J^2)=J\adj(J)$ by 
Lipman's version of the Brian{\c c}on-Skoda theorem, the latter 
scheme provides a "Gorensteinfication"
of $\Proj R_A(J)$. In fact, it turns out to be the canonical 
model of the minimal
desingularization $X\longrightarrow \Proj R_A(J)$ in the sense of \cite[Definition 0-3-11]{KMM}
(see Proposition~\ref{cano}). 
We remark that it remains to be true in any dimension that 
if $\Proj R_A(J)$ has only Gorenstein rational
singularities, then $J^n$ is an adjoint ideal 
for all large $n$ (see Proposition~\ref{general}). 
However, the blowup $\Proj R_A(\adj(J))$ 
need not be Gorenstein any more (see Example \ref{ex}).

Finally, we will give some examples in order to illustrate our results. 
Recall that an ideal is called
simple if it is not a product of two proper ideals. Lipman and 
Watanabe observed in~\cite[Proposition (2.3)]{LW} that a simple ideal 
is an adjoint ideal if and only if it is of order one. 
This implies in particular that it is generated by a regular sequence.  We
consider ideals with two simple factors and give in 
Theorem \ref{2simple} a necessary and sufficient
condition for this kind of an ideal to be an adjoint ideal. 

Throughout this paper, for an ideal $I$ in a local ring $(A,\fkm)$,
we denote by $\overline{I}$ its integral closure.
For an $\fkm$-primary ideal $I$, $\ell_A(A/I)$ is the length of $A/I$ and
$e(I)$ is the multiplicity of $A$ with respect to $I$.
We say that an ideal $Q\subset I$ is a reduction of $I$ 
if $I^{n+1} = QI^n$ for some $n$. A reduction is called minimal 
if it is minimal with respect to inclusion.
The Rees algebra of $I$ is $R_A(I) = \oplus_{n \geq 0} I^n$.
For any set $S$, $\# S$ stands for the number of elements of $S$.


\section{Preliminaries}
\label{sec:preli}

In this section, we will fix some notation and recall some basic facts 
concerning complete ideals in two-dimensional regular local rings,
which can be extracted from~\cite{L1}--\cite{L5} and~\cite{ZS}.

\medskip

\noindent 
{\bf (2.1)} \quad
Let $(A,\fkm)$ be a two-dimensional regular local ring with the fraction field $K$.
Two-dimensional regular local rings with fraction field $K$ are called {\it points}.
Let $x \in \fkm \setminus \fkm^2$, and let $M$ be a maximal ideal of $A[\fkm/x]$
containing $x$. The localization $A[\fkm/x]_M$ is then called a {\it quadratic transform\/} of $A$.
Quadratic transforms are points. Moreover, any two points $B\subset C$ can be 
connected by a unique sequence of quadratic transforms 
$$B =: B_0 \subset B_1 \subset \ldots \subset B_n:=C.$$ If
$\fkm_B$ and $\fkm_C$ denote the maximal ideals of $B$ and $C$, respectively,
then always $\fkm_B\subset \fkm_C$, and the residue field extension $B/\fkm_B \to C/\fkm_C$ 
has finite degree, denoted by $[C : B]$.

\medskip

\noindent 
{\bf (2.2)} \quad
The {\it order of an element\/} $0 \ne a \in A$ is  $\ord_A(a) = \max \{ n | a \in \fkm^n \}$.
The order function yields a discrete valuation $\ord_A:K \setminus \{ 0 \} \to \Bbb Z$.
The {\it order of an ideal\/} $I$ is $\ord_A(I) = \max \{ n | I \subset \fkm^n \}$. Let
$B\supset A$ be a point. Take the sequence of quadratic transforms $A =: A_0 \subset 
A_1 \subset \dots\subset A_n := B$. For an $\fkm$-primary ideal $I$, the {\it transform\/} 
$I^B$ of $I$ at $B$ is defined
inductively by setting first $I^A=I$, and then $I^{A_{i+1}} =  x^{-\ord_{A_i}(I^{A_i})}I^{A_i}
A_{i+1}$ for $i=0,\ldots,n-1$ where $x$ denotes a generator of the principal ideal  
$\fkm_{A_i} A_{i+1}$. It is $\fkm$-primary (unless it is the unit ideal). Moreover,
if $I$ is complete, then so is $I^B$. We write $r_B(I) = \ord_B(I^B)$ for short. A point $B\supset A$
with $r_B(I)>0$ is called a {\it base point\/} of $I$. The {\it support\/} of $I$,
$\Supp I = \{ B \mid r_B(I) > 0 \}$, is known to be a finite set.  The family $\B (I) = 
(r_B(I))_{B \supset A}$ is called the {\it point basis\/} of $I$.  A complete 
$\fkm$-primary ideal is completely determined by its point basis i.e.~for any
complete $\fkm$-primary ideals $I$ and $J$ in $A$, we have $I=J$ 
if and only if $\B(I)=\B(J)$.  Also note that $\B(I) \geq \B(J)$ implies $I \subset J$.

\medskip
\noindent 
{\bf (2.3)} \quad
We can always "blow up the base points" in order to arrive to a
non--singular scheme $X \rightarrow \Spec A$ such that $IO_X$ is invertible. More precisely,  
there is a sequence of
morphisms $$ X = X_{n+1}  \overset{f_n}\longrightarrow X_n \overset{f_{n-1}}\longrightarrow \ldots 
\overset{f_0}\longrightarrow X_0 = \Spec A$$
obtained by blowing up closed points $x_i\in X_i$ ($0 \leq i \le n$) such that $\Supp I=
\{\Oh_{X_i,x_i} \mid 0 \leq i  \le n\}$.  One can then give the point basis of an $\fkm$-primary ideal
the following geometric interpretation. Let $E_0,\ldots,E_n$ denote the exceptional divisors of the
morphism $X\rightarrow \Spec A$ i.e.~for every $0\le i\le n$ $E_i$ is the strict transform of $E'_i=f_i^{-1}(x_i)$ in the 
morphism $f_n\cdots f_{i+1}$.  The lattice $\mathbb Z E_0\oplus \ldots \mathbb \oplus 
\mathbb Z E_n$ now has besides $(E_0,\ldots,E_n)$ another basis, namely 
$(E^*_0,\ldots,E^*_n)$ where $E_i^*=(f_n\cdots f_{i+1})^*E'_i$ denotes the total 
transform of $E'_i$ ($0\le i\le n$). If $D$ is the effective Cartier divisor on $X$ for which $I\Oh_X=\Oh_X(-D)$,
then the point basis of $I$ just gives the coordinates of $D$ with respect to this basis. In other
words, we now have $D=\sum_{i=0}^nr_{\Oh_{X_i,x_i}}(I)E^*_i$ (see~\cite[Lemma 1.18]{C}).

If $Y = \Proj R_A(I)$,  then this morphism induces a morphism $X \to Y$ which
turns out to be the {\it minimal desingularization\/} of $Y$ in the sense that
every other desingularization $Z\to Y$ factors through $X$ (use~\cite[Corollary (27.3)]{L1}).  
Note that if $I$ is complete, then by~\cite[Proposition (1.2)]{L1} $Y$ has always only rational singularities.

\medskip 
\noindent
{\bf (2.4)} \quad
An ideal $I$ is called {\it simple\/} if it is not a product of two 
proper ideals. By Zariski's famous factorization theorem every complete ideal can be uniquely expressed as a product
of simple complete ideals. There is one to one correspondence between 
simple complete $\fkm$-primary ideals and points containing $A$
(see~\cite[p.~391, (E)]{ZS}). Indeed, 
given a point $B\supset A$ the corresponding simple ideal $\fkp_B$ is the unique
ideal $I$ in $A$ whose transform $I^B=\fkm_B$. Moreover, if
$A =: A_0 \subset A_1 \subset \ldots \subset A_n:=B$ is the quadratic
sequence associated to $B$, then $\Supp I=\{A_0,\ldots,A_n\}$.
We often write $\B(I) = (r_0, r_1,\ldots,r_n)$,  where $r_i = r_{A_{i}}(I)$.
Note that the following {\it reciprocity formula\/} holds: $\ord_A(I)=[B:A]\ord_B(\fkm_A)$
(see~\cite[(1.6.1)]{L3}). Also note 
that there is further one to one correspondence between points and prime divisors 
of $A$ given by $B\mapsto \ord_B$.  Recall here that a valuation $v$ of $K$ is called 
a {\it prime divisor\/} of $A$ if its valuation ring dominates $A$ 
and the transcendence degree of the residue field of the valuation ring over $A/\fkm$ 
is one. 

\medskip
\noindent 
{\bf (2.5)} \quad
A fundamental tool in this article is the notion of proximity.
Let $B$ and $C$ be points such that $C$ strictly contains $B$. 
We say that $C$ is {\it proximate\/} to $B$, and write $B \prec C$, 
if  $C$ is contained in the valuation ring associated to $\ord_B$.
We note that $B \prec C$ whenever $C$ is a quadratic transform of $B$.
More generally, in any chain $B =: B_0 \subset B_1 \subset \ldots \subset B_n$
of quadratic transformations $\{ B_j \mid B \prec B_j \} = \{ B_1, B_2,\ldots, B_m\}$
where $m$ is determined by the expression
$$\ord_{B_n}(\fkm_B) = (m-1) \ord_{B_n}(\fkm_{B_1}) + b \quad (0 < b \leq \ord_{B_n}(\fkm_{B_1}))$$
(see~\cite[Lemma (5.2.1)]{L3}).

\medskip
\noindent 
{\bf (2.6)} \quad
Let $0\neq (r_B)_{B \supset A}$ be a family of nonnegative integers,
with $r_B = 0$ for all but finitely many $B$.
There exists a unique complete $\fkm$-primary ideal $I$ in $A$ such that 
$r_B(I) = r_B$ for all $B \supset A$ if and only if the $proximity \ inequality$
$r_B \geq  \sum_{B \prec C} [C : B] r_C$ holds for 
each $B \supset A$ (see~\cite[Theorem 2.1]{L4}). The {\it excess\/} of the ideal $I$ at a point $B\supset A$ 
is defined as 
$$\Exc_B(I) := r_B(I) - \sum_{B \prec C}[C : B] r_C(I).$$
The proximity inequality implies that the excess is always nonnegative.
The unique expression of $I$ as a product of simple complete $\fkm$-primary 
ideals is then given by $$I = \prod_{\Exc_B(I)>0} \fkp_B^{\Exc_B(I)}$$ (see~\cite[Corollary 3.1]{L4}). 
\medskip

\noindent 
{\bf (2.7)} \quad
Let $f\colon X\to \Spec A$ be a proper birational morphism. If $X$ is normal, then the {\it relative 
canonical sheaf\/} $\omega_X$ can be defined as the dual of the relative Jacobian
sheaf ${\cal J}_X$ (see~\cite[p.~206, (2.3)]{LS}). We always think of $\omega_X$ as a subsheaf of the constant sheaf $K$. 
Note that in the case 
$X = \Proj S$ for some graded $A$-algebra $S$, one can think $\omega_X$ as the sheafification of the 
graded canonical module of $S$ (see, e.g.~, \cite[2.6.2]{HS}).
Let $K_X$ denote the {\it canonical divisor\/} of $X$, i.e.,  the Cartier 
divisor on $X$ for which $ \Oh_X(K_X)=\omega_X$. 

\medskip
\noindent 
{\bf (2.8)} \quad
Let $I$ be an ideal in $A$. Let $X \to \Spec A$ be a proper birational morphism such that $X$ is 
nonsingular and $I \Oh_X$ is invertible. Let $D$ be the effective divisor on $X$ satisfying 
$\Oh_X(-D)=I\Oh_X$. The module of global sections $\Gamma(X, I \omega_X)=\Gamma(X,\Oh_X(K_X-D))$ 
is then independent of the chosen desingularization $X \to \Spec A$. It is called the {\it adjoint\/} of $I$, 
and denoted by $\adj(I)$. This is clearly the same as the multiplier ideal of $I$ with the multiplying parameter
$c=1$ as defined in~\cite[Definition 9.2.3]{Laz}. One should note, however, that the term "adjoint"  is used in
a different sense in~\cite{Laz} (see~\cite[Definition 9.3.47]{Laz}).

It is known that $\adj(I)$ is a complete ideal in $A$ containing $I$. Evidently
also $\adj(I)=\adj(\overline I)$. 
Moreover, if $I$ is $\fkm$-primary, then so is $\adj(I)$ unless it is the unit ideal. 
The following fact (\cite[Proposition (3.1.2)]{L5}) will be essential for us in sequel:

\smallskip
{\it The point basis of $\adj(I)$ is given by $(\max \{ r_B(I)-1, 0 \})_{B\supset A}$.}

\begin{proof}
For the convenience of the reader we sketch a proof.  
For all points $B\supset A$, set 
\begin{eqnarray}
s_B &=& \left\{ 
\begin{array}{ll} r_B(I)-1 & \ B \in \Supp I; \\
0 & \ $otherwise.$ \\
\end{array} \right. \nonumber 
\end{eqnarray}
The numbers $s_B$ satisfy the proximity inequalities
\begin{align*}
s_B-\displaystyle \sum_{B \prec C}[C : B]s_C&\ge s_B-\displaystyle \sum_{{B \prec C} \atop {C \in \Supp I}}[C : B]s_C\\
&=\Exc_B(I)+\displaystyle \sum_{{B \prec C} \atop {C \in \Supp I}}[C : B]\ge 0.\end{align*}
By (2.6) there exists a complete ideal $J'\subset A$ such that $r_B(J')=s_B$ for all $B$.
Take a desingularization $X\rightarrow \Spec A$ as in (2.3) such that $I\Oh_X$ is invertible. Since $\Supp J'
\subset \Supp I$, $J'\Oh_X$ is invertible, too.  Moreover, if  $F$ is the effective divisor on $X$ with $\Oh_X(-F)=J'\Oh_X$, then
$F=(r_0-1)E^*_0+\ldots +(r_n-1)E^*_n$ where $r_i=r_{\Oh_{X,x_i}}(I)$ $(0\leq i \leq n)$.  Now $D=r_0E^*_0+\ldots +r_nE^*_n$. 
As $K_{X_{i+1}}=f_i^*K_{X_i}+E'_i$ for every $0\le
 i\le n$, we have $K_X=E^*_0+\ldots +E^*_n$.  So $D-K_X=F$. Since $J'$ is complete, this implies that $J'=\adj(I)$. 
\end{proof}

We will call an ideal an {\it adjoint ideal\/} if it is the adjoint of some ideal. According to Lipman's version of
the Brian{\c c}on--Skoda theorem~\cite[(2.3) and Conjecture (1.6)]{L5}
we have $\adj(I^n)=I^{n-1}\adj(I)$ for all $n\ge 1$.  A key role in the following 
will be played by the fact that in the case $X \to \Proj R_A(I)$ is the minimal desingularization, $I \omega_X$ 
is globally generated by $\adj(I)$, in other words, we have $I\omega_X=\adj(I)\Oh_X$~\cite[(3.1.1)]{L5}. The adjoint can also be 
calculated by means of a minimal reduction $Q\subset I$ using the formula $\adj(I) = Q:I$~\cite[Proposition (3.3)]{L5}. This
implies in particular that $\ell_A(A/\adj(I)) = \ell_A(I/Q)$.


\section{Main results}

Let $(A, \fkm)$ be a regular local ring of dimension two. We want to give criteria for a complete ideal 
$J$ in $A$ to be an adjoint. Evidently every principal ideal is an adjoint. Otherwise, the following
lemma shows that we can restrict ourselves to $\fkm$-primary ideals.

\begin{lemma} \label{firstlem}Let $J$ be a complete ideal in $A$ which is not principal. Write $J=xJ'$ for some 
element $x\in A$ and some $\fkm$-primary ideal $J'\subset A$. Then $J=\adj(I)$ for some ideal $I$ in $A$ if 
and only if $J'=\adj(I')$ for some $\fkm$-primary ideal $I'$ in $A$.  In particular, if an $\fkm$-primary ideal
is an adjoint ideal, then it is an adjoint of an $\fkm$-primary ideal. 
\end{lemma}

\begin{proof} By \cite[(1.2.3) (c)]{L5} we know that $\adj(xI)=x \adj(I)$ for any element $x\in A$ and any 
ideal $I\subset A$. If $J'=\adj(I')$ for some ideal $I'\subset A$, then this immediately implies that $J=\adj(I)$ 
for the ideal $I=xI'$. To prove the converse, we can write $I=yI'$ for some $y\in A$ and some $\fkm$-primary ideal 
$I'\subset A$. Then $xJ'=y\adj(I')$. Since the ideal $\adj(I')$ can't  now be a unit ideal, it is necessarily 
$\fkm$-primary. Hence $(x)=(y)$ so that $J'=\adj(I')$ as wanted.
\end{proof}

It is convenient to give conditions for a power of an ideal to be an adjoint:

\begin{lemma} \label{mainlemma}
Let $J$ be a complete $\fkm$-primary ideal in $A$. Let $n$ be a positive integer. The following conditions 
are then equivalent for a complete $\fkm$-primary ideal $I$ in $A$
\begin{enumerate}
\item[$(1)$] $J^n = \adj(I)$ and $\Supp I=\Supp J$;
\item[$(2)$] $J^{n+1} = I  \adj (J)$;
\item[$(3)$]  $r_B(I)=\begin{cases}
n \cdot r_B(J)+1 &B \in \Supp J; \cr
0                        &\hbox{\rm otherwise}.
\end{cases}$
\end{enumerate}
Moreover, if these equivalent conditions hold, then necessarily $I =  J^{n+1} : \adj(J)$.

\end{lemma}

\begin{proof} Recall from (2.8) that $r_B(\adj(I)) = \max \{ r_B(I) -1, 0 \}$ for all points $B\supset A$. 
The equivalence of $(1)$ and $(3)$ is then obvious. On the other hand, an ideal $I$ satisfies 
$J^{n+1} = I \cdot \adj (J)$ if and only if $r_B(I)=(n+1)r_B(J)-r_B(\adj(J))$ for all $B$, showing
the equivalence of $(2)$ and $(3)$. Conditions (1)--(3) are thus equivalent.  Since
$(J^{n+1} : \adj(J))\adj(J)=I \adj (J)$ for any complete ideal $I$ satisfying 
$J^{n+1} = I  \adj (J)$, the last statement follows from the completeness of $J^{n+1} : \adj(J)$.
\end{proof}

We can now prove

\begin{theorem}{\label{main}}
Let $J$ be a complete $\fkm$-primary ideal in $A$. Let $n$ be a positive integer. 
Then the following conditions are equivalent:
\begin{enumerate}
\item[$(1)$] We have $J^n = \adj(I)$ for some ideal $I$ in $A$;
\item[$(2)$] $\adj (J) \mid J^{n+1}$;
\item[$(3)$] $n \cdot \Exc_{B} (J) + 1 \geq 
    \displaystyle \sum_{{B \prec C} \atop {C \in \Supp J}}[C : B]$ 
for all points $B\in \Supp J$. 
\end{enumerate}
Moreover, if these equivalent conditions hold, then $I =  J^{n+1} : \adj(J)$ is the unique complete
$\fkm$-primary  ideal such that $J^n = \adj(I)$ and
$\Supp I=\Supp J$. 

\end{theorem}

\begin{proof} 
For all points $B\supset A$, set 
\begin{eqnarray}
s_B &=& \left\{ 
\begin{array}{ll}
n \cdot r_B(J)+1 & \ B \in \Supp J; \\
0 & \ $otherwise.$ \\
\end{array} \right. \nonumber 
\end{eqnarray}
It is now easily checked that  
$$s_B - \sum_{B \prec C}[C : B] s_C=n \cdot \Exc_{B} (J) + 1-  \displaystyle \sum_{{B \prec C} \atop {C \in \Supp J}}[C : B]$$ 
when $B\in \Supp J$. 
Note that by the proximity inequality $s_B=0$ implies that $s_C=0$ for all points $C$
proximate to $B$. 
But as explained in (2.6) the inequalities $s_B \geq \sum_{B \prec C}[C : B] s_C$ hold true for
every point $B$ if and only if there exists a complete $\fkm$-primary ideal $I$ in $A$ such that 
$r_B(I) = s_B$ for all $B$. In light of Lemma~\ref{mainlemma} it remains to prove that if $J^n = \adj(I')$ 
for some ideal $I'$ in $A$, then there exists a complete $\fkm$-primary ideal $I$ with
$\Supp I=\Supp J$ such that $J^n = \adj(I)$.  By Lemma~\ref{firstlem}  $I'$ is necessarily $\fkm$-primary. 
The formula $r_B(\adj(I')) = \max \{ r_B(I') -1, 0 \}$
implies that $r_B(I') = nr_B(J) + 1 = s_B$ for all $B \in \Supp J$. The proximity inequality now shows that 
$$s_B - \sum_{B \prec C}[C : B] s_C \ge r_B(I') - \sum_{B \prec C}[C:B]r_C(I') \geq 0$$
if $B \in \Supp J$. Hence there exists a complete ideal $I$ in $A$ with $r_B(I) = s_B$ for all $B$. 
The theorem has thus been proven.
\end{proof}

\begin{remark}{\label{rem2}} One can give also a more geometric condition for an ideal $J$ in $A$
to be an adjoint.  Indeed, let $f:X \to \Proj R_A(J)$ be the minimal desingularization of $\Proj R_A(J)$. 
Let $F$ be the effective divisor on $X$ such that $J \Oh_X = \Oh_X(-F)$. For any $B\in \Supp J$,
let $E_B$ denote the exceptional divisor on $X$ corresponding to $B$. Then the following
condition is equivalent to conditions (1)--(3) of Theorem~\ref{main} : 

\bigskip

\noindent
\begin{enumerate}
\item[$(4)$] $E_B^2\ge -2[B:A]+n(F\cdot E_B)$ for all points $B\in \Supp J$.
\end{enumerate}

\bigskip

\noindent
We sketch here a proof for the equivalence of (1) and (4). Suppose
that $J^n = \adj(I)$ for some 
complete $\fkm$-primary ideal $I$ in $A$.  
Recall from (2.8) that $\adj(I)\Oh_X=I\omega_X$.
Hence $J^n \Oh_X = I \omega_X$. Conversely, if this equation holds
for some complete $\fkm$-primary ideal $I$ in $A$ ,
then taking global sections gives $J^n = \adj(I)$. 
Now $J^n \Oh_X = I \omega_X$ means the same as $\Oh_X(-nF-K_X) = I \Oh_X$.
By the one to one correspondence between  complete $\fkm$-primary ideals in $A$ and 
anti-nef divisors on $X$ (see \cite[Section 18, p.~238]{L1}) this is further equivalent to $nF+K_X$ 
being anti-nef, the latter meaning that $(nF+K_X)\cdot E_B \le 0$ for all $B\in \Supp B$. 
But $K_X\cdot E_B=-2[B:A]-E_B^2$  by~\cite[Proposition (4.5.1)]{L3}. The above claim
thus follows. 
\end{remark}

\begin{corollary}{\label{power}}
Let $J$ be a complete $\fkm$-primary ideal in $A$.
If $J^n$ is an adjoint ideal for some positive integer $n$, then so is $J^k$ for any integer $k \geq \min(n, \ord_A(J)-1)$.
\end{corollary}

\begin{proof}
It is readily clear from condition (3) of Theorem \ref{main} that $J^k$ is an adjoint ideal
for all $k\ge n$.  To complete the proof, it is enough to prove that $J^{\ord_A(J)-1}$ is an adjoint ideal.  
We thus need to show that 
$$(\ord_A(J)-1) \cdot \Exc_B(J) + 1\geq \sum_{{B \prec C} \atop {C \in \Supp J}}[C : B]$$
for all $B \supset A$.

If $\Exc_B(J) = 0$, then this is already clear by the assumption. Suppose that $\Exc_B(J) > 0$.
Then $$(\ord_A(J)-1) \cdot \Exc_B(J) + 1 \geq \ord_A(J).$$
By the proximity inequality we now obtain
$$\ord_A(J) \geq r_B(J) \geq \sum_{B \prec C}[C:B] r_C(J) \geq 
\sum_{{B \prec C} \atop {C \in \Supp J}}[C : B]$$
as wanted.
\end{proof}

Note the following "dual version" of Lipman's Brian{\c c}on-Skoda theorem:

\begin{corollary}{\label{quotients}}
Let $J$ be a complete $\fkm$-primary ideal in $A$.
If $J^n$ is an adjoint ideal for some positive integer $n$, then $J^{n+p+1}:\adj(J)=J^p(J^{n+1}:\adj(J))$ for all non-negative integers $p$.
\end{corollary}

\begin{proof} Set $I_p= J^{n+p+1}:\adj(J)$ for all $p\ge 0$.  By Corollary~\ref{power} also $J^{n+p}$ is
an adjoint. Theorem~\ref{main} (2) therefore gives $J^{n+p+1}=I_p\adj(J)$. Similarly, $J^{n+1}=I_0\adj(J)$ so
that $I_p\adj(J)=J^pI_0\adj(J)$. Because all the ideals involved are complete, we get $I_p=J^p I_0$ as
wanted.
\end{proof}

Given an ideal $I$ in $A$, we define for every positive integer $r$ the $r$-th {\it iterated adjoint\/} $\adj^r(I)$ by setting
$$\adj^r(I)=\adj^{r-1}(\adj(I)).$$

Our purpose is to show next that high enough powers of an adjoint are always iterated adjoints. To this purpose, we
first need some lemmas. 

\begin{lemma}{\label{iteratedrfoldlem}} Let $I$ be a complete $\fkm$-primary ideal in $A$. Let $r$ be a positive integer. Let $f\colon X\to \Spec A$ be
a desingularization such that $I\Oh_X$ is invertible. Then $\adj^r(I)=\Gamma(X,I\omega_X^r)$.
\end{lemma}

\begin{proof} Set $\adj_r(I)=\Gamma(X,I\omega_X^r)$. Standard arguments imply that the definition
of this "$r$-fold adjoint"  is independent of
the choice of the desingularization $f$. In particular, we can assume that the induced morphism $X\to \Proj R_A(I)$ 
is the minimal desingularization of $\Proj R_A(I)$. But then $I\omega_X=\adj(I)\Oh_X$ by (2.8), and so $\adj_r(I)=\adj_{r-1}(\adj(I))$. 
Because $\adj^r(I)=\adj^{r-1}(\adj(I))$, the claim follows by induction.
\end{proof}

\begin{lemma}{\label{iteratedpowerprop}} Let $I$ be a complete $\fkm$-primary ideal in $A$. Then $\adj^r(I^r)=\adj(I)^r$ for every positive integer $r$.
\end{lemma}

\begin{proof} Let $f\colon X\to \Proj R_A(I)$ be the minimal desingularization. By (2.8) we know that $ I \omega_X
=\adj(I) \Oh_X$. Since $\adj(I)^r$ is complete and complete ideals are contracted (see~\cite[Proposition (6.2)]{L1}), 
Lemma~\ref{iteratedrfoldlem} then gives 
$\adj^r(I^r)=\Gamma(X,I^r\omega_X^r)=\Gamma(X,\adj(I)^r\Oh_X)
=\adj(I)^r$.
\end{proof}

We are now ready to prove the following: 

\begin{proposition}{\label{Fittingprop}} Let $J$ be a complete $\fkm$-primary ideal in $A$. Let $r$ be a positive integer. If $J$ is an adjoint, then 
$J^n$ is an $r$-th iterated adjoint for every $n\ge r$. In particular, $J^n=F_{r+1}(I)$ for some ideal $I$ in $A$.
\end{proposition}

\begin{proof} Suppose that $J=\adj(I)$ for some ideal $I\subset A$. By Lemma~\ref{firstlem} $I$ is necessarily $\fkm$-primary.
Moreover, as $\adj(I)=\adj(\overline I)$, we can assume that $I$ is complete. 
Since $\adj^n(I^n)=\adj^r(\adj^{n-r}(I^n))$ for $n\ge r$, the first claim
is now a consequence of Lemma~\ref{iteratedpowerprop} according to which $\adj^r(I^n)=\adj(I)^n$. 
The last statement is then~\cite[Proposition (3.16) b) ]{HuSw}.
\end{proof}

Next we will connect adjoint ideals to the Gorenstein property of blowups. We first need the following lemma which
describes the canonical sheaf of a blowup along a complete ideal.

\begin{lemma}{\label{Iwy}}
Let $I$ be a complete $\fkm$-primary ideal in $A$.
Set $Y = \Proj R_A(I)$.
Then $I \omega_Y = \adj(I) \Oh_Y$.
\end{lemma}

\begin{proof}
Let $f: X \to Y $ be the minimal desingularization. Since $Y$ has only rational singularities (see (2.3)), we have $f_*(\omega_X)=\omega_Y$.
By the projection formula we now get $f_*(I \omega_X) = I\omega_Y$. On the other hand, an application of
Lemma~\ref{pushrfoldlem} below in the case $r=1$ then yields $I\omega_Y=\adj(I) \Oh_Y$ proving the claim.
\end{proof}

\begin{lemma}{\label{pushrfoldlem}} Let $I$ be a complete $\fkm$-primary ideal in $A$. Let $r$ be a positive integer. Let $f\colon X\to Y$ be
the minimal desingularization of $Y=\Proj R_A(I)$. Then $f_*(I^r\omega_X^r)=\adj(I)^r\Oh_Y$.
\end{lemma}

\begin{proof} Recalling again from (2.8) that $ I \omega_X=\adj(I) \Oh_X$, we get $f_*(I^r \omega_X^r) =f_*(\adj(I)^r\Oh_X)$.
As $Y$ has only rational singularities (see (2.3)),  the completeness of $\adj(I)$ implies 
by~\cite[Proposition (6.5)]{L1} that of $\adj(I) \Oh_Y$. By~\cite[Theorem (7.1)]{L1} $\adj(I)^r \Oh_Y$ is then complete, too. 
But complete ideals are contracted (see~\cite[Proposition (6.2)]{L1}) so that 
$f_* (\adj(I)^r \Oh_X)=\adj(I)^r \Oh_Y$ as wanted.
\end{proof}

\begin{theorem}{\label{gor}}
Let $J$ be a complete $\fkm$-primary ideal in $A$, and let $n$ be a positive integer.
Set $Y=\Proj R_A(J)$. Then the following conditions are equivalent.
\begin{enumerate}
\item[$(1)$] $Y$ is Gorenstein and $J^n \omega_Y^{-1}$ is 
generated by global sections;
\item[$(2)$] $J^n$ is an adjoint ideal.
\end{enumerate}
\end{theorem}

\begin{proof} The implication (1) $\Rightarrow$ (2) holds also in higher dimensions, and is a special case of 
Proposition~\ref{general} below. Extra assumptions made there are of course not needed in dimension 
two. Note that by (2.3) $Y$ has only rational singularities.
It remains to prove the implication (2) $\Rightarrow$ (1). By Theorem \ref{main} we can find a complete 
$\fkm$-primary
ideal $I$ in $A$ satisfying $J^{n+1} = I \adj(J)$. From Lemma \ref{Iwy} we get $J \omega_Y = \adj(J) \Oh_Y$. 
Hence $J^{n+1} \Oh_Y = 
I \Oh_Y \cdot \adj(J) \Oh_Y = I \Oh_Y \cdot J \omega_Y$. 
Therefore $\omega_Y$ is invertible (using~\cite[Theorem 11.6 d]{E}). Since $Y$ has only rational singularities, it is 
also Cohen-Macaulay. Hence $Y$ is Gorenstein.
Moreover, we see that $J^n \omega_Y^{-1} = I \Oh_Y$.
As $I$ is complete, this implies that $I=\Gamma(Y,J^n \omega_Y^{-1})$.
Hence $J^n \omega_Y^{-1}$ is globally generated.
\end{proof}

\begin{proposition}{\label{general}} Let $A$ be an excellent regular local ring of equicharacteristic zero, 
and let $J$ be a normal ideal in $A$.  If $Y = \Proj \R_A(J)$ is a Gorenstein scheme with 
only rational singularities and the sheaf $J^n\omega_Y^{-1}$ is generated by global sections, then $J^n$ 
is an adjoint ideal. In particular, this holds for all large enough $n$.  
\end{proposition}

\begin{proof} Let $f\colon X \to Y$ be a desingularization. The inclusions $A=\Gamma(X,\omega_X) \subset 
\Gamma(Y,\omega_Y)\subset A$ imply that $\Gamma(Y,\omega_Y)=A$. Then $\Oh_Y \subset \omega_Y$ 
and $\omega_Y^{-1}\subset \Oh_Y$.  As $J^n\omega_Y^{-1}$ is generated by global 
sections, we thus have $J^n \omega_Y^{-1} =I \Oh_Y$ for the ideal $I:=\Gamma(Y,J^n\omega_Y^{-1})\subset A$.
Note that $I \Oh_Y$ and $I \Oh_X$ are invertible so that 
$I \Oh_X = f^*(I \Oh_Y)$. Since $Y$ has only rational singularities, we have $f_*\omega_X= \omega_Y$. The projection 
formula now gives $$\adj(I) = \Gamma (X, I \omega_X) = \Gamma (Y, I f_*(\omega_X))=\Gamma(Y,I\omega_Y)=\Gamma(Y,J^n\Oh_Y)=J^n.$$
In the last step we used the assumption that $J^n$ is integrally closed (see, e.g. , \cite[(1.4)]{LT}).
\end{proof}

The implication (2) $\Rightarrow$ (1) in Theorem~\ref{gor} does not hold true in general.  This comes out from
the following example suggested to us by T.~J\"arvilehto.

\begin{example}{\label{ex}}
{\rm
Set $A = k[x,y,z]_{(x,y,z)}$ where $k[x,y,z]$ is a polynomial ring 
over a field $k$. Consider the ideal $I = \overline{(x,y^3,z^3)}$.
Then one can check that $I$ and $I^2$ are integrally closed.
By \cite[Proposition 3.1]{RRV}, 
it follows that $I^n$ is integrally closed for all $n$.
In particular, since $I$ is monomial, $\Proj \R_A(I)$ then has only 
rational singularities.
By a formula due to Howald (\cite[Main Theorem]{H}), 
$J := \adj(I^2)$ is generated by monomials $x^a y^b z^c$ with
$a + \frac{b}{3} + \frac{c}{3} > \frac{1}{3}$, whence 
$J = \overline{(x,y^2,z^2)}$.
On the other hand, $\R_A(J) = A[xt,y^2t, yzt,z^2t]$, where $t$
is an indeterminate. It is easy to check that the coordinate ring $A[y^2/x, yz/x,z^2/x]$
of the affine chart $D_+(xt) \subset \Proj \R_A(J)$ is not Gorenstein.
So $J$ is an adjoint ideal, but $\Proj \R_A(J)$ is not
Gorenstein.
We note that $\Proj \R_A(J)$ has only rational singularities by the same argument as above, since
$J$ and $J^2$ are integrally closed.
}
\end{example}

Let $(A,\fkm)$ be a two-dimensional regular local ring.
It is well known that if $I$ is a complete $\fkm$-primary
ideal in $A$, then the blowup $Y=\Proj R_A(I)$ is Gorenstein if and only if its singularities are
double points (for an algebraic argument, see~\cite[Corollary 1.6]{HuSa}). Classification of rational double points in terms of
"configuration diagrams" can be found in~\cite[p.~258]{L1}. Recall here that $Y$ has always rational singularities (see (2.3)).
In fact, being a so called sandwiched singularity each singularity of $Y$ is 
necessarily of type $(A_r)$ for some $r$ (see~\cite[p.~426]{S}).  What we want to do here is to characterize the Gorenstein 
property of $Y$ in terms
of adjoints.

\begin{corollary}{\label{gor_ord}}
Let $J$ be a complete $\fkm$-primary ideal.
Then the following are equivalent.
\begin{enumerate}
\item[$(1)$] $\Proj \R_A(J)$ is Gorenstein;
\item[$(2)$] $J^n$ is an adjoint ideal for some positive integer $n$ (and then for all integers $n\ge \ord_A(J)-1$).
\end{enumerate}
\end{corollary}

\begin{proof} 
Set $Y = \Proj \R_A(J)$.
Because $J\Oh_Y=\Oh_Y(1)$ is ample, 
$J^n \omega_Y^{-1}$ is generated by global sections for large $n$. 
Using Theorem \ref{gor} we then know that $\Proj \R_A(J)$ is Gorenstein if
and only if $J^n$ is an adjoint ideal for large $n$. But by
Corollary \ref{power} the latter is equivalent to $J^n$
being an adjoint ideal for some (and then for all) $n\ge \ord_A(J)-1$.
\end{proof}

The following proposition is now an analogy of Theorem~\ref{main}:

\begin{proposition}{\label{primdiv}}
Let $J$ be a complete $\fkm$-primary ideal in $A$.
Then the following conditions are equivalent.
\begin{enumerate}
\item[$(1)$] $\Proj \R_A(J)$ is Gorenstein;
\item[$(2)$] $\adj(J)\mid J^n$ for some positive integer $n$ (and then for all integers $n\ge \ord_A(J)$);
\item[$(3)$] For any point $B\in \Supp B$ with $\Exc_{B} (I) = 0$, there is a
unique point $C\in \Supp B$ proximate to $B$. Moreover, then $[C:B]=1$. 
\end{enumerate}
\end{proposition}

\begin{proof} Because of Corollary~\ref{gor_ord}, the equivalence of (1) and (2) is immediately clear from
condition (2) of Theorem~\ref{main}. It remains to prove the equivalence of (2) and (3). We use
condition (3) of Theorem~\ref{main}.  This clearly holds for the finitely many points $B\in \Supp J$ with 
$\Exc_{B} (J) >0$ when $n$ is large. It is thus only necessary to consider points $B\in \Supp B$ with 
$\Exc_{B} (J) = 0$. Observe that by the definition of excess there then exists some $C\in \Supp J$ proximate to $B$.
Our condition now reads 
$$1 =   \displaystyle \sum_{{B \prec C} \atop {C \in \Supp J}}[C : B],$$
but this is clearly equivalent to (3).
\end{proof}

\begin{remark} Let us state some geometric variants of condition (3), which perhaps look more familiar.  
Set $Y=\Proj R_A(J)$. Let $f\colon X
\to Y$ be the minimal desingularization.  For any point $B\in \Supp J$, let $E_B$ denote the exceptional 
divisor on $X$ corresponding to $B$. The divisors $E_B$ with $\Exc_B(J)>0$ corresponding to the irreducible
components of the closed fiber of the morphism $Y\to \Spec A$,  we observe that $\Exc_{B} (J) = 0$ if and 
only $f$ contracts $E_B$ to a closed point on $Y$ i.e.~$E_B$
is $f$-exceptional. An application of the formula $$E_B^2=[B:A](-1-\sum_{{B \prec C} \atop {C \in \Supp J}}[C : B])$$
(see~\cite[p.~232]{L3})
then shows 
that an equivalent version of condition (3) is 
\medskip

\noindent
(3') $E_B^2=-2[B:A]$ for every point $B\in \Supp J$ such that $E_B$ is $f$-exceptional. 
\medskip

\noindent
The formula $K_X\cdot E_B=-E_B^2-2[B:A]$ (see~\cite[Proposition (4.5.1)]{L3}) gives a further variant:

\medskip
\noindent
(3'') $K_X\cdot E_B=0$ for every point $B\in \Supp J$ such that $E_B$ is $f$-exceptional.

\end{remark}

Corollary~\ref{gor_ord} immediately gives the following:

\begin{corollary}{\label{cor}}
Let $I$ be a complete $\fkm$-primary ideal in $A$. 
Then $\Proj \R_A(\adj(I))$ is Gorenstein.
\end{corollary}

Let $I$ be a complete $\fkm$-primary ideal in a two-dimensional 
regular local ring $(A,\fkm)$.  Set $Y = \Proj \R_A(I)$.
By Lipman's version of the Brian{\c c}on-Skoda theorem (see (2.8)), 
we have $\adj(I^2) = I \adj(I)$.
Combining this to Corollary \ref{cor} yields that $Z = \Proj \R(I \adj(I))$ 
is Gorenstein.
By the universal property of blowing up
$Z$ coincides with the blowup $\PProj \oplus_{n \geq 0} \adj(I)^n \Oh_Y$
of $Y$ along $\adj(I)\Oh_Y$.
This implies that $Z$ can be considered as a 
"Gorensteinfication" of $Y$. Let us connect this observation to the notion
of  "canonical model". Recall that in dimension two rational Gorenstein
singularities are exactly the canonical singularities. 
Let $f:X \to Y$ be the minimal desingularization.
The sheaf of $\Oh_Y$-algebras $\Rees_{X/Y} = \oplus_{n \geq 0}
f_*(\omega_X^{\otimes n})$ is called the {\it canonical algebra} of $X$ 
over $Y$. The corresponding scheme $\PProj \Rees_{X/Y}$
is the {\it canonical model\/} of $X$ over $Y$ (\cite[Definition 0-3-11]{KMM}).

\begin{lemma}{\label{canlem}}
We have $Z\cong \PProj \Rees_{X/Y}$ as $Y$-schemes.
In particular, $\PProj \Rees_{X/Y}$ is Gorenstein.
\end{lemma}

\begin{proof} Lemma~\ref{pushrfoldlem} combined with the projection formula implies that 
$$f_*(\omega_X^{\otimes n}) \otimes_{\Oh_Y} (I \Oh_Y)^n =f_*(I^n\omega_X^{\otimes n}) = \adj(I)^n \Oh_Y.$$  
Therefore
$$\displaystyle \PProj \Rees_{X/Y} \cong 
\PProj \bigoplus_{n \geq 0}( f_*(\omega_X^{\otimes n}) 
\otimes_{\Oh_Y} (I \Oh_Y)^n) \cong \PProj \bigoplus_{n \geq 0} \adj(I)^n \Oh_Y.$$
But then
$$\PProj \Rees_{X/Y} \cong \Proj \R_A(I \adj(I)) = \Proj \R_A(\adj(I^2)).$$
\end{proof}

The following proposition is now a two-dimensional 
local analogy of \cite[Theorem 0-3-12]{KMM}.

\begin{proposition}{\label{cano}}
The following conditions are equivalent.
\begin{enumerate}
\item[$(1)$] $Y = \Proj \R_A(I)$ is a Gorenstein scheme;
\item[$(2)$] The structure morphism $\pi : \PProj \Rees_{X/Y} \to Y$ is an isomorphism.
\end{enumerate}
\end{proposition}

\begin{proof}
Suppose that $Y$ is Gorenstein. Then $\adj(I) \Oh_Y$ is invertible by Lemma \ref{Iwy}
so that $Z\cong \PProj \oplus_{n \geq 0} \adj(I)^n \Oh_Y \cong Y$. By Lemma \ref{canlem} there is
an isomorphism $Z\cong \PProj \Rees_{X/Y}$ of $Y$-schemes. Therefore $\pi$ is an isomorphism.  
Conversely, if $\pi$ is an isomorphism, then $Y$ is Gorenstein by Lemma \ref{canlem}.
\end{proof}


\section{Examples}
\label{sec:ex}

In this section we will give several examples of adjoint ideals. For simplicity, we
mainly consider ideals having at most two simple factors. We continue to assume that 
$(A, \fkm)$ is a two-dimensional regular local ring. The following proposition is a slight 
improvement of 
\cite[Proposition 2.3]{LW}.

\begin{proposition}{\label{L_W}}
Let $I$ be a simple complete $\fkm$-primary ideal in $A$. 
Then the following conditions are equivalent:
\begin{enumerate}
\item[$(1)$] $I^n$ is an adjoint ideal for all positive integers $n$;
\item[$(2)$] $I^n$ is an adjoint ideal for some positive integer $n$;
\item[$(3)$] $\ord_A(I) = 1$.
\end{enumerate}
\end{proposition}

\begin{proof}
(1) $\Rightarrow$ (2): Trivial.

\noindent
(2) $\Rightarrow$ (3): Suppose that $I^n$ is an adjoint ideal. By Theorem \ref{main} (2)
$\adj(I) H = I^{n+1}$ for some ideal $H\subset A$ . Since $I$ is simple, this implies
that $\adj(I) = I^q$ for some $q \geq 0$. Recall from (2.8) that $r_B(\adj(I))=\max\{r_B(I)-1,0\}$
for all points $B\supset A$. Therefore we must have $\ord_A(I) - 1 = q \cdot \ord_A(I)$
so that $\ord_A(I) = 1$.

\noindent
(3) $\Rightarrow$ (1): By Corollary \ref{power} 
it is enough to check that $I$ is an adjoint ideal. Indeed,
if $\ord_A(I) = 1$, then $I = \adj(I^2)$, since $\B(I) = (1,1,\ldots,1)=\B(\adj(I^2))$.

\end{proof}

We next consider ideals having two simple factors. We first make the following basic observation:

\begin{lemma}{\label{ordI=1}}
Suppose that $A/\fkm$ is an algebraically closed field.
Let $I$ and $J$ be simple complete $\fkm$-primary ideals in $A$.
Let $p$ and $q$ be positive integers.
If $I^pJ^q$ is an adjoint ideal, then
either $\Supp I \subset \Supp J$ or $\Supp I \supset \Supp J$.
\end{lemma}

\begin{proof}
The condition that $I^pJ^q$ is an adjoint implies by Theorem \ref{main} (3) that
$$\Exc_B(I^pJ^q) + 1 \geq 
\# \{ C \in \Supp IJ \mid B \prec C \}$$ 
for all $B \supset A$, because $A/\fkm$ is algebraically closed.
Let $A=:A_0\subset A_1\subset \ldots \subset A_s$ and 
$A=:B_0\subset B_1\subset \ldots \subset B_t$ denote the base points of
$I$ and $J$, respectively.  Suppose, for example, that $s \leq t$. 
If $s=0$, then $I = \fkm$, whence $\Supp I = \{ A \}$, and we are done.
Let $s>0$.
Then $A_1$, $B_1 \in \{ C \in \Supp IJ \mid A \prec C \}$.
Since $\Exc_A(I^pJ^q) = 0$, this is impossible unless $A_1 = B_1$.
Proceeding inductively, we finally conclude that
$A_i =  B_i$ for $i = 1,2,\ldots,s$, which implies that 
$\Supp I \subset \Supp J$.
\end{proof}

We are now able to characterize those adjoint ideals which have two simple factors:

\begin{theorem}{\label{2simple}}
Suppose that $A/\fkm$ is an algebraically closed field.
Let $I \ne J$ be simple complete $\fkm$-primary ideals  in $A$. 
Let $p$ be a positive integer. Set $k = \ord_A(J)$.
Then the following conditions are equivalent:
\begin{enumerate}
\item[$(1)$] $I^pJ^q$ is an adjoint ideal for some positive integer $q$ and $\Supp I \subset \Supp J$;
\item[$(2)$] $I^pJ^q$ is an adjoint ideal for all positive integers $q$ and $\Supp I \subset \Supp J$; 
\item[$(3)$] 
	\begin{enumerate}
	\item[$(i)$] $p+1 \geq k$,
	\item[$(ii)$] $\Supp I \subset \Supp J$,
            \item[$(iii)$] $\B(J) = (k,\ldots,k,1,\ldots,1)$ where $k$ appears 
$\# \Supp I$ times. 
	\end{enumerate}
Moreover, these conditions imply $\B(I) = (1,1,\ldots,1)$;
\item[$(4)$]
	\begin{enumerate}
	\item[$(i)$] $p+1 \geq k$,
        \item[$(ii)$] $\ord_A(I) = 1$,
	\item[$(iii)$] $\adj(J) = I^{k-1}$;
	\end{enumerate}
\item[$(5)$]
	\begin{enumerate}
	\item[$(i)$] $p+1 \geq k$,
        \item[$(ii)$] $\ord_A(I) = 1$,
	\item[$(iii)$] $J \subset I^k$,
	\item[$(iv)$] $e(J)-\ell_A(A/J) = \ell_A(A/I^{k-1})$.
	\end{enumerate}
\end{enumerate}
\end{theorem}

\begin{proof}
(1) $\Leftrightarrow$ (3) and (2) $\Leftrightarrow$ (3):
We can assume that $\Supp I \subset \Supp J$.
Let $A=:A_0 \subset \ldots \subset A_t$ denote the base points of $J$. Let  
$\Supp I = \{ A_0, A_1,\ldots,A_s \}$.  As $I\not=J$, $s<t$.
By Theorem \ref{main} (3) $I^pJ^q$ is an adjoint ideal if and only if 
$$\# \{ A_j  \mid A_i \prec A_j \}\le  \Exc_{A_i}(I^pJ^q) +1$$ for
all $i=0,\ldots,t$. This is the same as
$$\# \{ A_j  \mid A_i \prec A_j \}\le \begin{cases}1\ &\hbox{if $s\not=i<t$;}\cr p+1\ &\hbox{if $s=i$.}\end{cases}\qquad (\ast)$$
Recall that by reciprocity (see (2.4))
$r_{A_i}(J)=\ord_{A_t}(\fkm_{A_i})$ for all $i=0,\ldots,t$.  It now follows from (2.5) that $(\ast)$ is equivalent to having
$$k=r_{A_0}(J)=\ldots=r_{A_s}(J)\quad\hbox{and}\quad  r_{A_{s+1}}(J)=\ldots =r_{A_t}(J)=1$$ where
$p+1\ge k$. The desired equivalences thus follow.  
Moreover, if these equivalent conditions hold, then an application of (2.5) again
gives $r_{A_0}(I)=\ldots=r_{A_s}(I)=1$. 

\noindent
(3) $\Leftrightarrow$ (4): Note that $\adj(J)= I^{k-1}$ implies $\Supp I\subset \Supp J$, since always $J\subset \adj(J)$. 
In any case we can thus assume that $\Supp I\subset \Supp J$. Moreover,  $\ord_A(I)=1$ means the same as $\B(I)=
(1,\ldots,1)$. Then $\adj(J)= I^{k-1}$ if and only if $\B(\adj(J))=(k-1,\ldots,k-1)$. 
But recalling from (2.8) that $r_B(\adj(J))=\max\{r_B(J)-1,0\}$ for all points $B\supset A$, the latter is equivalent to 
$\B(J) = (k,\ldots,k,1,\ldots,1)$.
Thus (3) and (4) are equivalent.

\noindent
(4) $\Leftrightarrow$ (5):  
Since $\B(I)=(1,\ldots,1)$, we observe that $\adj(J)=I^{k-1}$ implies $J\subset I^k$. Conversely, if
$J\subset I^k$, then $\adj(J)\subset \adj(I^k)=I^{k-1}$. The equality $\adj(J)=I^{k-1}$ holds
if and only if $\ell_A(A/\adj(J))=\ell_A(A/I^{k-1})$. Now recall from (2.8) that $\ell_A(A/\adj(J))=\ell(J/Q)$
for any minimal reduction $Q\subset J$. As $e(J)=\ell(A/Q)$, the equivalence of (4) and (5) follows.

\end{proof}

Theorem~\ref{2simple} implies in particular that if a complete $\fkm$-primary ideal having two simple factors
is an adjoint, then one factor always has order one. It then follows that it is of the form $(x,y^r)$
for a suitable regular system of parameters $\{x,y\}$ of $A$. The next proposition now gives information 
about the other factor.

\begin{proposition}{\label{e(J)}}
Suppose that $A/\fkm$ is an algebraically closed field. 
Let $I$ and $J$ be simple complete $\fkm$-primary ideals  in $A$ with $\ord_A(I)=1$ and $k=\ord_A(J) > 1$. 
Let $p$ and $q$ be positive integers. Set $s = \# \Supp I$.  If
$I^pJ^q$ is an adjoint ideal, then $e(J) \geq k^2s+k$. Moreover,
$e(J) = k^2s+k$ if and only if $J = \overline{(x^k, y^{ks+1})}$ 
for a suitable regular system of parameters $\{ x,y \}$ of $A$. 
\end{proposition}

\begin{proof} Theorem~\ref{2simple} (3)  implies that we now have $\Supp I\subset \Supp J$.
Let $A=:A_0 \subset \ldots \subset A_t$ denote the base points of $J$, where
$\Supp I=\{A_0,\ldots,A_{s-1}\}$. We then get moreover, 
$\B(I)=(1,\ldots,1)$ and $\B(J)=(k,\ldots,k,1,\ldots,1)$ where $k$ appears
$s$ times. Using (2.5) we see that there must be $k$ points $A_s,\ldots,A_{s+k-1}$
proximate to $A_{s-1}$. Therefore $t-s+1\ge k$. It now follows
from the formula of Hoskin and Deligne (see \cite[Corollary (3.8)]{L2}) that
$e(J) = k^2s+t-s+1 \geq k^2s+k$.

Suppose then that $e(J)=k^2s+k$. Hence $t=s+k-1$. 
Because $\ord_A(I)=1$, we know that $I=(x,y^r)$ for some regular system of parameters
$\{x,y\}$ of $A$. By the uniqueness of the quadratic sequence corresponding to
$I$, we then necessarily have $r=s$ and $$A_i=A_{i-1}[x/y^i]_{(y,x/y^i)}\qquad (i=
1,\dots,s-1).$$ Set $v=\ord_{A_t}$. Since $x/y\in \fkm_{A_1}$, we have $v(x)>v(y)$. 
As $v(\fkm)= \ord_A(J)=k$ by reciprocity (see (2.4)), this implies that $v(y)=k$.  Furthermore, we conclude 
from $x/y^{s-1}\in \fkm_{A_{s-1}}$ that $v(x/y^{s-1})\ge v(y)$. But then $A_s=A_{s-1}[x/y^s]_{\fkn}$ where $\fkn$
denotes the center of $v$ on $A_{s-1}[x/y^s]$. As $A/\fkm$ is algebraically closed, $\fkn=
(y,x/y^s-\lambda)$ for some $\lambda \in A\setminus \fkm$. Set $x'=x-\lambda y^s$.
Then $\{x',y\}$ is a system of parameters of $A$ with $I=(x',y^s)$. By changing the
notation, we can assume that $x=x'$. Because $v(\fkm_{A_s})=r_{A_s}(J)=1$, we must have 
$v(x/y^s)=1$ i.e.~$v(x)=ks+1$. It then follows that $$A_i=A_{i-1}[y/(x/y^s)^{i-s}]_{(x/y^s,y/
(x/y^s)^{i-s})}\qquad (i=s+1,\dots,s+k-1).$$ It is now easily checked that if $J'=\overline{(x^k,y^{sk+1})}$,
then $\B(J')=\B(J)$ implying that $J=J'$. Conversely, if $J=J'$, then by the above $t=s+k-1$ so that
$e(J)=k^2+k$. 
\end{proof}

The next corollary is an analogy of Proposition \ref{L_W}.

\begin{corollary}{\label{ordJ=1}}
Suppose that $A/\fkm$ is an algebraically closed field.
Let $I$ and $J$ be simple complete $\fkm$-primary ideals  in $A$ with 
$\ord_A(I) = \ord_A(J)= 1$. Let $p$ and $q$ be positive integers.
Then the following conditions are equivalent:
\begin{enumerate}
\item[$(1)$] $(I^pJ^q)^n$ is an adjoint ideal for all positive integers $n$;
\item[$(2)$] $(I^pJ^q)^n$ is an adjoint ideal for some positive integer $n$; 
\item[$(3)$] Either $I \subset J$ or $I \supset J$ holds true.
\end{enumerate}
\end{corollary}

\begin{proof}
(1) $\Rightarrow$ (2): Trivial.

\noindent
(2) $\Rightarrow$ (3): By Lemma \ref{ordI=1} 
$\Supp I \subset \Supp J$ or $\Supp I \supset \Supp J$.
Since $\ord_A(I) = \ord_A(J)= 1$, this is equivalent to saying that $\B(I) \geq \B(J)$ or 
$\B(I) \leq \B(J)$. Hence $I \subset J$ or $I \supset J$.

\noindent
(3) $\Rightarrow$ (1): Suppose, for example, that $\Supp I\subset \Supp J$. 
The claim is then a consequence of Theorem \ref{2simple} (3), because $np+1 \geq 1$ for all $n >0$.
\end{proof}

A complete $\fkm$-primary ideal $I$ in $A$ is of order one if and only if the adjoint $\adj(I)=A$. 
We want next to consider the class of complete ideals for which 
$\adj(I)=\fkm$. To this purpose, we recall the general notion of minimal multiplicity of an ideal introduced by Goto in~\cite{G2}. Let $I$ be 
an $\fkm$-primary ideal in a Cohen-Macaulay local ring $(A, \fkm)$ which has an infinite residue field. Let $\mu_A(I)$ stand for the least 
number of generators of $I$. One then always has the inequality $e(I)\ge \mu_A(I)+\ell(A/I)-\dim A$. If the equality holds, then $I$ is said to
have {\it minimal multiplicity}. This condition is equivalent to saying that $\fkm I \subset Q$ for any minimal reduction $Q$ of $I$ 
(\cite[Lemma (2.1)]{G2}). We now return to the two-dimensional situation.

\begin{proposition}{\label{minimal}}
Suppose that $A/\fkm$ is an infinite field.
Let $I$ be a complete $\fkm$-primary ideal in $A$. 
Then the following conditions are equivalent:
\begin{enumerate}
\item[$(1)$] $I$ has minimal multiplicity;
\item[$(2)$] $\adj(I)=A$ or $\adj(I)=\fkm$;
\item[$(3)$] One of the following conditions holds true.
\begin{enumerate}
\item[$(i)$] $\ord_A(I) = 1$,
\item[$(ii)$] $I$ is simple and $\B(I) = (2,1,\ldots,1)$,
\item[$(iii)$] $I = JH$ with $\ord_A(J)=\ord_A(H)=1$ and 
$\Supp J \cap \Supp H = \{ A \}$.
\end{enumerate}
\end{enumerate}

\end{proposition}

\begin{proof}

(1) $\Leftrightarrow$ (2): As mentioned above $I$ has minimal multiplicity if and only if
$\fkm I \subset Q$ for any minimal reduction $Q\subset I$ i.e.~$\fkm \subset Q:I$. But 
$\adj(I)=Q:I$ (see 2.8), and  we are done. 

(2) $\Rightarrow$ (3): Recall that $r_B(\adj(I)) = \max \{ r_B(I)-1,0 \}$ for all points 
$B\supset A$. So $\adj(I) = A$ implies $\ord_A(I)=1$. In the case $\adj(I) =
\fkm$ it follows that $r_A(I) = 2$ and $r_B(I) \leq 1$ for all $B \ne A$. If $I$ is simple, then
$\Supp I$ forms a chain, and we get $\B(I) = (2,1,\ldots,1)$. Otherwise we can write $I=JH$ 
for some proper ideals $J$ and $H$ in $A$. Then both $J$ and $H$ are necessarily of 
order one. Moreover, $r_B(J) + r_B(H) =r_B(I) \leq 1$ for any $B\not=A$. Hence 
$\Supp J \cap \Supp H = \{A\}$.

(3) $\Rightarrow$ (2):
The formula $r_B(\adj(I)) = \max \{ r_B(I)-1,0 \}$ shows that in the case (i)
$\adj(I)=A$ whereas in the case (ii) or (iii) we have $\Supp \adj(I)=\{A\}$ and 
$\ord_A(\adj(I))=1$ so that $\adj(I)=\fkm$.
\end{proof}

\begin{corollary}{\label{min_adj}}
Suppose that $A/\fkm$ is an infinite field.
Let $I$ be a complete $\fkm$-primary ideal in $A$. 
If $I$ has minimal multiplicity, 
then the following conditions are equivalent:
\begin{enumerate}
\item[$(1)$] $I$ is an adjoint ideal;
\item[$(2)$] Either $\ord_A(I) =1$ or $I = \fkm J$ for some ideal $J$.
\end{enumerate}
\end{corollary}

\begin{proof}
\noindent
(1) $\Rightarrow$ (2):
If $I$ is simple, then by Proposition \ref{L_W} it follows  that
$\ord_A(I) =1$.
Suppose that $I$ is not simple.
Then by Proposition \ref{minimal} $I = JH$ where $J$ and $H$ are simple of order one and
$\Supp J \cap \Supp H = \{ A \}$.
On the other hand, by Lemma \ref{ordI=1}, 
$\Supp J \subset \Supp H$ or $\Supp J \supset \Supp H$, 
which implies that $J = \fkm$ or $H = \fkm$.

\noindent
(2) $\Rightarrow$ (1):
If $\ord_A(I) =1$, then $I$ is a simple ideal, whence
$I$ is an adjoint ideal by Proposition \ref{L_W}.
Suppose that $I = \fkm J$ for some ideal $J\subset A$.
By Proposition \ref{minimal} $\ord_A(J) = 1$.
Hence $I$ is an adjoint ideal by Corollary \ref{ordJ=1}.
\end{proof}

\begin{corollary}{\label{sim_min_adj}}
Suppose that $A/\fkm$ is an infinite field.
Let $I$ be a simple complete $\fkm$-primary ideal in $A$. 
Then the following conditions are equivalent:
\begin{enumerate}
\item[$(1)$] $I$ has minimal multiplicity;
\item[$(2)$] $\ord_A(I) \leq 2$ and 
$\fkm^p I^q$ is an adjoint ideal for some (and then for all) positive integers $p$ and $q$.
\end{enumerate}
\end{corollary}

\begin{proof}
\noindent
(1) $\Rightarrow$ (2):
By Corollary \ref{ordJ=1} we may assume that $\ord_A(I) > 1$ .
Then $\B(I) = (2,1,\ldots,1)$ by Proposition \ref{minimal}. In particular,
we have $\ord_A(I) = 2$. It then follows from Theorem \ref{2simple}
(3) that $\fkm^p I^q$ is an adjoint ideal for any $p,q > 0$ .

\noindent
(2) $\Rightarrow$ (1):
If $I$ is an ideal of order one, then $I$ has minimal multiplicity.
When $\ord_A(I) = 2$, Theorem \ref{2simple} (3) implies that $\B(I) = (2,1,\ldots,1)$. 
According to Proposition \ref{minimal} $I$ has thus minimal multiplicity.
\end{proof}

\bigskip

\noindent
{\bf Acknowledgements:}
We are grateful to Tarmo J\"arvilehto for useful discussions. The first
and third author were supported by the Academy of Finland, project
48556. The second author would like to thank the University of 
Helsinki for the hospitality during the preparation of this article.

\end{document}